\documentclass
[12pt]{article} \oddsidemargin 0in \evensidemargin 0in \textheight
9in \textwidth 6.5in \topmargin 0in \headheight 0in
\usepackage{amsmath,amsfonts,amsthm,amssymb}
\usepackage{hyperref}

\newtheorem{thm}{Theorem}
\newtheorem{theorem}{Theorem}[section]
\newtheorem{cor}[theorem]{Corollary}

\newtheorem{defn}{Definition}

\title{ Remarks on Multiplicative Metric Spaces and Related Fixed Points}

\author {K. Abodayeh$^1$\and A. Pitea$^2$ \and W. Shatanawi$^3$\thanks{(Permanent address) Hashemite University, Zarqa-Jordan,
Email:swasfi@hu.edu.jo},\and T. Abdeljawad$^4$
   \\ {\it $^{1,3,4}$Department of Mathematics and General Sciences}
   \\{ Prince Sultan University}
   \\ {\it Riyadh, Saudi Arabia 11586 }
   \\ {\it E-mails: kamal@psu.edu.sa, wshatanawi@psu.edu.sa, tabdejawad@psu.edu.sa}\\{\it Faculty of Applied Sciences}
   \\{\it $^2$Department of Mathematics and Informatics}
   \\{University Politehnica of Bucharest}
 \\{\it Bucharest, Romania}
 \\{\it E-mail:arianapitea@yahoo.com}}

\date{}

\begin{document}
 \maketitle{}

\begin{abstract}

In this article we studied the relationship between metric spaces
and multiplicative metric spaces. Also, we pointed out some fixed
and common fixed point results under some contractive conditions in
multiplicative metric spaces can be obtained from the corresponding
results in standard metric spaces.
\end{abstract}

\section{Introduction}
The notion of multiplicative metric space was introduced by Bashirov
et al. \cite{bashirov}. In 2012, Ozavsar and Cervikel \cite{ozavsar}
defined the notion of convergence in multiplicative metric spaces
and studied some fixed point results in such space. After that, many
 researchers consider this space and many results on fixed point theory were considered. \\

We begin by introducing the definition of multiplicative metric
space.

\begin{defn}\cite{bashirov} Let $X$ be a nonempty set and $p: X \times X \longrightarrow
[1,+\infty).$ We say that $(X,p)$ is a multiplicative metric space
if for all $ x,y,z \in X$ we have:
\begin{enumerate}
\item $p(x,y)\geq 1$ and $x=y$ if and only if $p(x,y)=1;$

\item $p(x,y)= p(y,x);$

\item $p(x,z)\le p(x,y)p(y,z).$
\end{enumerate}
\end{defn}

 The definition of a multiplicative Cauchy sequence in
 a multiplicative metric space is given as follows:
\begin{defn}\cite{ozavsar}
A sequence $\{ x_n \}$ in  a multiplicative metric space $(X,p)$ is
said to be multiplicative Cauchy sequence if for all $\epsilon >1$,
there exists $N\in{\cal N}$ such that $p(x_n ,x_m )<\epsilon$ for
all $m,n\geq N$. Also, if every multiplicative Cauchy sequence is
convergent, then $(X,p)$ is called a complete multiplicative metric
space.
\end{defn}
For the definitions of open balls and convergence in multiplicative
metric spaces, we refer the reader to \cite{ozavsar}.\\

Ozavsar and Cervikel \cite{ozavsar} introduced the concept of
multiplicative contraction and proved that every multiplicative
contraction in a complete multiplicative metric space has a unique
fixed point.

\begin{defn}\label{d2}\cite{ozavsar}
Let $(X,p)$ be a multiplicative metric space. A mapping
$f:X\rightarrow X$ is called {\it multiplicative contraction} if
there exists a real number $\lambda\in [0,1)$ such that $p(f(x_1
),f(x_2 ))\leq p(x_1 ,x_2 )^{\lambda}$ for all $x_1 ,x_2 \in X$.
\end{defn}

\begin{thm}\label{t1}\cite{ozavsar}\label{banach}
Let $(X,p)$ be a complete multiplicative metric space and let
$f:X\rightarrow X$ be a multiplicative contraction. Then $f$ has a
unique fixed point.
\end{thm}

The notion of weakly commuting mappings was introduced by Sessa
\cite{sessa} in 1982. While, Jungck \cite{jungck}
 initiated the concept of weakly compatible mappings in 1996 as a generalization of the notion of weakly commuting mappings.  \\
Moreover, many authors studied many fixed point theorems for weakly
commuting mappings in metric spaces. See \cite{jungck}, \cite{sessa1}, \cite{sessa2}, \cite{agrawal}, \cite{moosaei}, and \cite{kumam}.

\begin{defn}\cite{jungck}
Let $A$ and $S$ be self-mappings on a metric space $( X, d)$. Then,
$A$ and $S$ are said to be weakly compatible if they commute at
their coincident point; that is, $Ax =Sx$ for some $x\in X$ implies
$ASx= SAx$.
\end{defn}

\begin{defn}\cite{sessa}
Let $S$ and $T$ be two self-mappings of a metric space $(X,d)$. Then
$S$ and $T$ are said to be weak commutative mappings if
$$
d(STx,TSx)\leq d(Sx,Tx),
$$
for all $x\in X$.
\end{defn}
It is clear that if $S$ and $T$ are weak commutative mappings, then
$S$ and $T$ are weakly compatible.

He et al. \cite{xiaoju} employed the concept of weakly commutative
mappings to introduce and prove the following common fixed point
theorem in multiplicative metric spaces.

\begin{thm}\label{t2}\cite{xiaoju}
Let $(X, p)$ be a complete multiplicative metric space. Suppose that
$A, B, S$ and $T$ are four self-mappings of $X$ satisfying the
following conditions:
\begin{enumerate}
\item $T(X)\subseteq A(X)$ and  $S( X)\subseteq B( X)$;
\item The pairs $(S, A)$ and $(T, B)$ are weakly commutative;
\item One of $A, B, S$ and $T$ is continuous;
\item
\begin{eqnarray}\label{e1}
p(Sx,Ty)\leq \{ \max \{ p(Ax,By), p(Ax,Sx), p(By, Ty), p(Ax,Ty),
p(By,Sx) \} \}^{\lambda}, \lambda\in\Large (0,\frac{1}{2}\Large ).
\end{eqnarray}
\end{enumerate}
Then A, B, S and T have a unique common fixed point.
\end{thm}

In this paper, we study the relationship between the multiplicative
metric space and the standard metric space. Also, we show that the
proof of Theorem \ref{t1} and Theorem \ref{t2} are obtained from the
corresponding results in standard metric spaces.

\section{Main Results}

We start by giving the relationship between the multiplicative
metric space and the standard metric space.\\

If we have a multiplicative metric space $(X,p)$, then the
corresponding metric space $(X,d_p )$ is given by the following
theorem.

\begin{thm}
Let $(X,\rho)$ be a multiplicative metric space. Define $d_p:X\times
X\rightarrow [0,+\infty)$ by
\begin{eqnarray*}
d_p (x,y) =\ln (p(x,y)).
\end{eqnarray*}
Then $(X,d_p)$ is a metric space.
\end{thm}
{\bf{Proof.}} Follows from the properties of
logarithms.\indent$\Box$\\

 Moreover, if we have a
metric space $(X,d)$, then the corresponding multiplicative metric
space $(X,p_d )$ is given by the following theorem.

\begin{thm}
Let $(X,d)$ be a metric space. Define $p_d:X\times X\rightarrow
[0,+\infty)$ by
 $$p_d (x,y)=e^{d(x,y)}.$$
Then $(X,p_d)$ is a multiplicative metric space.
\end{thm}
{\bf{Proof.}} The proof follows from the properties of exponential
functions.\indent$\Box$\\

Now, we can transfer and  prove many applications considered on
multiplicative metric space to a standard metric space and use
their proof in the metric space case. \\
For instance, considering the definition of contraction on
multiplicative metric space, the proof the following result is a
straightforward.\\

\begin{thm}
A sequence  $\{ x_n \}$ is a multiplicative Cauchy sequence in  a
multiplicative metric space $(X,p)$ if and only if $\{ x_n \}$ is a
Cauchy sequence in the corresponding metric space $(X,d_p )$.
\end{thm}

Applying the logarithmic function to the multiplicative contraction
inequality that have been defined in Definition \ref{d2}, will give
us the inequality\\

$$d_p  (f(x_1 ),f(x_2 ))=\ln p(f(x_1 ),f(x_2 )) \leq \lambda d_p (x_1 ,x_2
).$$\\

Note that if $(X,p)$ is a complete multiplicative metric space, then the corresponding metric space $(X,d_p )$ is
also a complete metric space. \\

It is obvious we can get the regular contraction inequality which
was introduce by Banach. Therefore, one can prove
the result in Theorem \ref{banach} using the new metric space $(X,d_p )$ and the Banach contraction theorem. \\
\\
 We have furnished all the necessary backgrounds to
present the proof of Theorem \ref{t1} from the standard metric space.\\
\\
 \noindent{\bf Proof of Theorem \ref{banach}:}\\

Since $(X,p)$ is a complete multiplicative metric space, the corresponding metric space $(X,d_p )$ is a complete metric space.
 Also, since $f$ is a multiplicative contraction, it is a contraction in $(X,d_p )$. Therefore it satisfies the
 Banach contraction conditions and thus it has a unique fixed
 point.\indent$\Box$\\

Now, we furnish all the necessary backgrounds to prove Theorem
\ref{t2} from the corresponding standard metric space.

Recall the following definition.
\begin{defn}\cite{7}
Let X be a nonempty set and let $d : X \times X \longrightarrow
[0,+\infty)$ be a function satisfying the following conditions:
\begin{enumerate}
\item $d ( x, y) = d ( y, x)$.
\item If $d (x, y) =  0$ then we have $x = y$.
\item  $d (x, y) \leq d (x, z) + d (z, y)$ for all $x, y, z\in X$.
\end{enumerate}
Then the pair $( X ,d )$ is called the d-metric space. It also appeared under the name of metric-like space \cite{Harandi}.
\end{defn}

It is clear that every metric space is a d-metric space.

\begin{defn}\cite{7}
A sequence $\{ x_n \}$ in a d-metric space $( X ,d )$ is called a
Cauchy sequence if for given
$\epsilon > 0$, there exists $N \in {\cal N}$ such that $d(x_n ,x_m )<\epsilon$ for all  $m, n \geq N$.\\
Also, $( X, d )$ is called complete if every Cauchy sequence in it
is convergent.
\end{defn}

\begin{defn}
A function $\phi :[0,\infty ) \rightarrow [0,\infty )$ is said to be
contractive modulus if $\phi (t)<t$ for $t>0$.
\end{defn}

\begin{defn}
A real valued function $\phi$ defined on $X \subset {\cal R}$ is
said to be upper semicontinuous if
$$
\lim_{n\rightarrow\infty}\phi (t_n )\leq \phi (t),
$$
for every sequence $\{t_n \}$ with $\lim_{n\rightarrow\infty} t_n =
t$.
\end{defn}

Panthi et al.\cite{panthi} introduced and proved the following
result.

\begin{thm}\label{t3}\cite{panthi}
Let $(X, d)$ be a complete d-metric space. Suppose that $A, B, S$
and $T$ are four self mappings of $X$ satisfying the following
conditions:
\begin{enumerate}
\item $T(X)\subseteq A(X)$ and  $S( X)\subseteq B( X)$,
\item $d(Sx,Ty)\leq \phi (m(x,y))$ where $\phi$ is an upper semicontinuous contractive modulus and
$$
m(x,y)\leq \max \{ d(Ax,By), d(Ax,Sx), d(By, Ty),
\frac{1}{2}d(Ax,Ty), \frac{1}{2}d(By,Sx) \}
$$
\item The pairs $(S, A)$ and $(T, B)$ are weakly compatible.
\end{enumerate}
Then A, B, S and T have a unique common fixed point.
\end{thm}

Now, we utilize Theorem \ref{t3} to introduce and prove the
following result.
\begin{cor}\label{main}
Let $(X, d)$ be a complete metric space. Suppose that $A, B, S$ and
$T$ are four self mappings of $X$ satisfying the following
conditions:
\begin{enumerate}
\item $T(X)\subseteq A(X)$ and  $S( X)\subseteq B( X)$,
\item suppose there exists $\lambda\in[0,1)$ such that
 $$d(Sx,Ty)\leq \lambda \max \{ d(Ax,By), d(Ax,Sx), d(By, Ty), \frac{1}{2}d(Ax,Ty), \frac{1}{2}d(By,Sx)
 \}$$
\item The pairs $(S, A)$ and $(T, B)$ are weakly compatible.
\end{enumerate}
Then A, B, S and T have a unique common fixed point.
\end{cor}

{\bf{Proof.}} Define $\phi : {\cal R}\longrightarrow {\cal R}$ by
$$\phi(t)=\lambda t,$$
where $0<\lambda <1$. Then
\begin{enumerate}
\item $\phi$ is a contractive modulus, and
\item $\phi$ is an upper semicontinuous on ${\cal R}$.
\end{enumerate}
From Theorem \ref{t3}, we get the result.\indent$\Box$\\

Now, we are ready to present the proof of Theorem \ref{t2} from the
corresponding standard metric space.\\

\noindent{\bf Proof of Theorem \ref{t2}:}\\

Take $\ln$ to both side in Inequality \ref{e1} in Theorem \ref{t2},
we get
$$d_p(Sx,Ty)\leq
\lambda\max\{d_p(Ax,By),d_p(Ax,Sx),d_p(By,Ty),d_p(Ax,Ty),d_p(By,Sx)\},$$
where $\lambda\in(0,\frac{1}{2})$ which is a special case of the
contractive condition in Corollary \ref{main}. Since every pair of
weakly commutative mappings is weakly compatible. All the hypotheses
of Corollary \ref{main} hold. Thus the four mappings $A,B,S$ and $T$
have a unique common fixed point.\indent$\Box$\\

\noindent{\bf {\large Conclusion:}}\\ \\
From our discussions, we note that some fixed and common fixed point
theorems in multiplicative metric spaces can be deduced from the
corresponding standard metric spaces. Moreover, we can formulate
many fixed and common fixed point theorems in multiplicative metric
spaces from the corresponding results in standard metric spaces. So
the researchers must be careful in working in multiplicative metric
spaces.


\begin{thebibliography}{999}
\bibitem{panthi} D. Panthi, K. Jha, P. Jha, and P. Kumari, ``A common fixed point theorem for two pairs of mappings in dislocated metric space,'' {\it American Journal of Computational Mathematics,} {\bf 5}, (2015) 106-112.

\bibitem{jungck} G. Jungck, ``Common fixed points for noncontinuous nonself maps on nonmetric spaces,'' {\it Far Eastern Journal of Mathematics and Science} {\bf 4(2)}(1996)199-215.

\bibitem{sessa}  S. Sessa,  ``On a weak commutativity condition of mappings in fixedpoint consideration,'' {\it Inst. Math.} {\bf 32} (1982) 149 - 153.

\bibitem{bashirov}  AE. Bashirov, EM. Kurplnara, A. Ozyaplcl, ``Multiplicative calculus and its applications,''
                    {\it Universal Journal of Mathematical Analysis and Applications} {\bf 337}
                  (2008), 36-48.
\bibitem{xiaoju}  X. He, M. Song, D. Chen , ``Common fixed points for weak commutative mappings on a multiplicative metric space,''
                    {\it  Fixed Point Theory and Applications } {\bf 2014:48}
                  (2014).

\bibitem{ozavsar}  M.~Ozavsar, and AC. Cervikel , ``Fixed points of multiplicative contraction mappings on multiplicative metric spaces,''
                    \url{http://arxiv.org/abs/1205.5131v1} (2012).

\bibitem{abbas}  M.~Abbas, B. Ali, and Y. Suleiman ``Common fixed points of locally contractive mappings in multiplicative metric spaces with applications,'' {\it International Journal of Mathematics and Mathematical Sciences} {\bf 2015}
                  (2015), Article ID 218683.

\bibitem{sarwar}  M.~Sarwar, and B. Rome , ``Some Unique fixed points theorems in  multiplicative metric space,''
                    \url{http://arxiv.org/abs/1410.3384} (2014).

 \bibitem{7}  P. Hitzler, and A.K. Seda, ``Dislocated Topologies,'' {\it Journal of Electrical Engineering}, {\bf 51} (2000) 3-7.

\bibitem{sessa1}  S. Sessa, ``On a weak commutativity condition of mappings in fixed point considerations,'' {\it Publ. Inst. Math.}, {\bf 32} (1982) 149-153.

\bibitem{sessa2}  S. Sessa and B. Fisher, ``On common fixed points of weakly commuting mappings and set-valued mappings,'' {\it International Journal of Mathematics and Mathematical Sciences}, {\bf 9} (1986) 323-329.

\bibitem{agrawal}  R.P. Agarwal, R.K. Bisht, and N. Shahzad, ``A comparison of various noncommuting conditions in metric fixed point theory and their applications,'' {\it Fixed Point Theory and Applications}, {\bf 2014} (2014) 2014:38.

\bibitem{moosaei}  M. Moosaei, ``Common fixed points for some generalized contraction pairs in convex metric spaces,'' {\it Fixed Point Theory and Applications}, {\bf 2014} (2014) 2014:98.

\bibitem{kumam}  P. Kumam, W. Sintunavarat, S. Sedghi, and N. Shobkolaei, ``Common Fixed Point of Two R -Weakly Commuting Mappings in b-Metric Spaces,'' {\it Journal of Function Spaces}, {\bf 2015} (2015) Article ID 350840.

\bibitem{Harandi}  A.A. Harandi, ``Metric-like spaces, partial metric spaces and fixed points,'' {\it Fixed Point Theory and Applications}, {\bf 2012} (2012) 2012:204.

\end{thebibliography}
\end{document}